\documentclass[a4paper]{article}
\pdfoutput=1

\usepackage{graphicx,wrapfig}
\usepackage{lipsum}
\usepackage{amsmath,amssymb,mathrsfs}
\usepackage{hyperref}
\usepackage{url}
\usepackage{mathtools}
\usepackage{changepage}
\usepackage{multirow}
\usepackage{wrapfig}
\usepackage{subfig}
\usepackage{color}
\usepackage{epsfig,enumerate,amsmath,amsfonts,amssymb,amsthm,mathrsfs,ifpdf}
\usepackage{indentfirst,relsize}
\usepackage{setspace,graphicx}
\usepackage{latexsym}
\usepackage[all]{xy}
\usepackage[usenames,dvipsnames]{pstricks}
\usepackage{pst-grad} 
\usepackage{pst-plot} 
\usepackage[margin = 2.50cm]{geometry}

\usepackage{algorithm}

\usepackage{algpseudocode}

\newcommand{\remove}[1]{}

\newtheorem{theo}{Theorem}
\newtheorem{lem}[theo]{Lemma}
\newtheorem{pre}[theo]{Proposition}

\newtheorem{cl}[theo]{Claim}

\newtheorem{remk}{Remark}

\title{Algorithmic study of $d_2$-transitivity of graphs}

\author{Subhabrata Paul \and Kamal Santra}

\author{Subhabrata Paul\footnote{Department of Mathematics, IIT Patna, India} \and Kamal Santra\footnote{Department of Mathematics, IIT Patna, India} }

\date{}

\begin{document}

\maketitle
\begin{abstract}
	Let $G=(V, E)$ be a graph where $V$ and $E$ are the vertex and edge sets, respectively. For two disjoint subsets $A$ and $B$ of $V$, we say $A$ \emph{dominates} $B$ if every vertex of $B$ is adjacent to at least one vertex of $A$. A vertex partition $\pi = \{V_1, V_2, \ldots, V_k\}$ of $G$ is called a \emph{transitive partition} of size $k$ if $V_i$ dominates $V_j$ for all $1\leq i<j\leq k$. In this article, we initiate the study of a generalization of transitive partition, namely \emph{$d_2$-transitive partition}. For two disjoint subsets $A$ and $B$ of $V$, we say $A$ \emph{$d_2$-dominates} $B$ if, for every vertex of $B$, there exists a vertex in $A$, such that the distance between them is at most two. A vertex partition $\pi = \{V_1, V_2, \ldots, V_k\}$ of $G$ is called a \emph{$d_2$-transitive partition} of size $k$ if $V_i$ $d_2$-dominates $V_j$ for all $1\leq i<j\leq k$. The maximum integer $k$ for which the above partition exists is called \emph{$d_2$-transitivity} of $G$, and it is denoted by $Tr_{d_2}(G)$. The \textsc{Maximum $d_2$-Transitivity Problem} is to find a $d_2$-transitive partition of a given graph with the maximum number of parts. We show that this problem can be solved in linear time for the complement of bipartite graphs and bipartite chain graphs. On the negative side, we prove that the decision version of the \textsc{Maximum $d_2$-Transitivity Problem} is NP-complete for split graphs, bipartite graphs, and star-convex bipartite graphs.

\end{abstract}


\section{Introduction}

Partitioning a graph is one of the fundamental problems in graph theory. In the partitioning problem, the objective is to partition the vertex set (or edge set) into some parts with desired properties, such as independence, minimal edges across partite sets, etc. In literature, partitioning the vertex set into certain parts so that the partite sets follow particular domination relations among themselves have been studied. Let $G$ be a graph with $V(G)$ as its vertex set and $E(G)$ as its edge set. When the context is clear, $V$ and $E$ are used instead of $V(G)$ and $E(G)$. The \emph{neighbourhood} of a vertex $v\in V$ in a graph $G=(V, E)$ is the set of all adjacent vertices of $v$ and is denoted by $N_G(v)$. The \emph{degree} of a vertex $v$ in $G$, denoted as $\deg_G(v)$, is the number of edges incident to $v$. A vertex $v$ is said to \emph{dominate} itself and all its neighbouring vertices. A \emph{dominating set} of $G=(V, E)$ is a subset of vertices $D$ such that every vertex $x\in V\setminus D$ has a neighbour $y\in D$, that is, $x$ is dominated by some vertex $y$ of $D$. For two disjoint subsets $A$ and $B$ of $V$, we say $A$ \emph{dominates} $B$ if every vertex of $B$ is adjacent to at least one vertex of $A$.

There has been a lot of research on graph partitioning problems that are based on a domination relationship between the different sets. Cockayne and Hedetniemi introduced the concept of \emph{domatic partition} of a graph $G=(V, E)$ in 1977, in which the vertex set is partitioned into $k$ parts, say $\pi =\{V_1, V_2, \ldots, V_k\}$, such that each $V_i$ is a dominating set of $G$ \cite{cockayne1977towards}. The number that represents the highest possible order of a domatic partition is referred to as the \emph{domatic number} of G, and it is denoted by $d(G)$. Another similar type of partitioning problem is the \emph{Grundy partition}. Christen and Selkow introduced a Grundy partition of a graph $G=(V, E)$ in 1979 \cite{CHRISTEN197949}. In the Grundy partitioning problem, the vertex set is partitioned into $k$ parts, say $\pi =\{V_1, V_2, \ldots, V_k\}$, such that each $V_i$ is an independent set and for all $1\leq i< j\leq k$, $V_i$ dominates $V_j$. The maximum order of such a partition is called the \emph{Grundy number} of $G$, and it is denoted by $\Gamma(G)$. In 2018, J. T. Hedetniemi and S. T. Hedetniemi \cite{hedetniemi2018transitivity} introduced a transitive partition as a generalization of the Grundy partition. A \emph{transitive partition} of size $k$ is defined as a partition of the vertex set into $k$ parts, say $\pi =\{V_1,V_2, \ldots, V_k\}$, such that for all $1\leq i< j\leq k$, $V_i$ dominates $V_j$. The maximum order of such a transitive partition is called \emph{transitivity} of $G$ and is denoted by $Tr(G)$. Recently, in 2020, Haynes et al. generalized the idea of domatic partition as well as transitive partition and introduced the concept of \emph{upper domatic partition} of a graph $G$, where the vertex set is partitioned into $k$ parts, say $\pi =\{V_1, V_2, \ldots, V_k\}$, such that for each $i, j$, with $1\leq i<j\leq k$, either $V_i$ dominates $V_j$ or $V_j$ dominates $V_i$ or both \cite{haynes2020upper}. The maximum order of such an upper domatic partition is called \emph{upper domatic number} of $G$, and it is denoted by $D(G)$. All these problems, domatic number \cite{chang1994domatic, zelinka1980domatically, zelinka1983k}, Grundy number \cite{effantin2017note, furedi2008inequalities, hedetniemi1982linear, zaker2005grundy, zaker2006results}, transitivity \cite{hedetniemi2018transitivity, haynes2019transitivity, paul2023transitivity, santra2023transitivity}, upper domatic number \cite{haynes2020upper, samuel2020new} have been extensively studied both from an algorithmic and structural point of view. Clearly, a Grundy partition is a transitive partition with the additional restriction that each partite set must be independent. In a transitive partition $\pi =\{V_1, V_2, \ldots, V_k\}$ of $G$, we have domination property in one direction, that is, $V_i$ dominates $V_j$ for $1\leq i< j\leq k$. However, in a upper domatic partition $\pi =\{V_1,V_2, \ldots, V_k\}$ of $G$, for all $1\leq i<j\leq k$, either $V_i$ dominates $V_j$ or $V_j$ dominates $V_i$ or both. The definition of each vertex partitioning problem ensures the following inequalities for any graph $G$. For any graph $G$, $1\leq \Gamma(G)\leq Tr(G)\leq D(G)\leq n$.

In this article, we introduce a similar graph partitioning problem, namely \emph{$d_2$-transitive partition}, which is a generalization of transitive partition. For two disjoint subsets $A$ and $B$, we say $A$ \emph{$d_2$-dominates} $B$ if, for every vertex of $B$, there exists a vertex in $A$, such that the distance between them is at most two. A \emph{$d_2$-transitive partition} of size $k$ is defined as a partition of the vertex set into $k$ parts, say $\pi =\{V_1,V_2, \ldots, V_k\}$, such that for all $1\leq i< j\leq k$, $V_i$ $d_2$-dominates $V_j$. The maximum order of such a $d_2$-transitive partition is called \emph{$d_2$-transitivity} of $G$ and is denoted by $Tr_{d_2}(G)$. The \textsc{Maximum $d_2$-Transitivity Problem} and its corresponding decision version are defined as follows:

\noindent\textsc{\underline{Maximum $d_2$-Transitivity Problem(M$d_2$TP)}}

\noindent\emph{Instance:} A graph $G=(V,E)$

\noindent\emph{Solution:} An $d_2$-transitive partition of $G$

\noindent\emph{Measure:} Order of the $d_2$ transitive partition of $G$\\

\noindent\textsc{\underline{Maximum $d_2$-Transitivity Decision Problem(M$d_2$TDP)}}

\noindent\emph{Instance:} A graph $G=(V,E)$, integer $k$

\noindent\emph{Question:} Does $G$ have a $d_2$-transitive partition of order at least $k$?\\

Note that every transitive partition is also a $d_2$-transitive partition. Therefore, for any graph $G$, $1\leq Tr(G)\leq Tr_{d_2}(G)\leq n$. Moreover, the difference between $Tr_{d_2}(G)$ and $Tr(G)$ can be arbitrarily large. For complete bipartite graphs $K_{t,t}$, $Tr_{d_2}(G)=2t$ whereas $Tr(G)=t+1$.
From the complexity point of view, there are some graph classes where transitivity can be solved in linear time, but $d_2$-transitivity is NP-complete. For example, in split graphs, the transitivity problem can be solved in linear time \cite{santra2023transitivity}, but later in this paper, we show that $d_2$-transitivity is NP-complete in split graphs. There are some vertex partition parameters where the value of the parameter in a subgraph can be greater than the original graph. The upper domatic number is one such example. But in the case of a $d_2$-transitive partition, $Tr_{d_2}(H)\leq Tr_{d_2}(G)$, for every subgraph $H$ of $G$. As a consequence, for a disconnected graph, the $d_2$-transitivity is equal to the maximum $d_2$-transitivity among all of its components. Therefore, we focus only on connected graphs in this paper.

In this paper, we study the computational complexity of the $d_2$-transitivity of graphs. The main contributions are summarized below:
\begin{enumerate}
	
	\item[1.] We show that the $d_2$-transitivity can be computed in linear time for the complement of bipartite graphs and bipartite chain graphs.\\
	
	\item [2.] We show that the \textsc{M$d_2$TDP} is NP-complete for split graphs, bipartite graphs, and star-convex bipartite graphs.
\end{enumerate}

The rest of the paper is organized as follows. Section 2 contains basic definitions and notations that are followed throughout the article. We have also discussed the basic properties of $d_2$-transitivity of graphs in this section. Section 3 describes linear-time algorithms for the complement of bipartite graphs and bipartite chain graphs. In Section 4, we show the \textsc{M$d_2$TDP} is NP-complete in split graphs, bipartite graphs, and star-convex bipartite graphs. Finally, Section 5 concludes the article.

\section{Preliminaries}

\subsection{Definitions and notations}
Let $G=(V, E)$ be a graph with $V$ and $E$ as its vertex and edge sets, respectively. A graph $H=(V', E')$ is said to be a \emph{subgraph} of a graph $G=(V, E)$ if and only if $V'\subseteq V$ and $E'\subseteq E$. For a subset $S\subseteq V$, the \emph{induced subgraph} on $S$ of $G$ is defined as the subgraph of $G$ whose vertex set is $S$ and edge set consists of all of the edges in $E$ that have both endpoints in $S$, and it is denoted by $G[S]$. The \emph{complement} of a graph $G=(V,E)$ is the graph $\overline{G}=(\overline{V}, \overline{E})$, such that $\overline{V}=V$ and $\overline{E}=\{uv| uv\notin E\}$. The \emph{open neighbourhood} of a vertex $x\in V$ is the set of vertices $y$ that are adjacent to $x$, and it is denoted by $N_G(x)$. The \emph{closed neighborhood} of a vertex $x\in V$, denoted as $N_G[x]$, is defined by $N_G[x]=N_G(x)\cup \{x\}$. For any $x, y\in V$, the \emph{distance} between $x$ and $y$ is defined as the number of edges in the shortest path starting at $x$ and ending at $y$ in $G$, and it is denoted by $d(x,y)$. The \emph{diameter} of a graph $G$ is defined as the greatest length of the shortest path between each pair of vertices, and it is denoted by $diam(G)$. Let $G$ be a graph; the \emph{square} of $G$ is a graph with the same set of vertices as $G$ and for any $x,y\in V$, $xy$ is an edge in the square graph if and only if $d(x,y)\leq2$. The square graph of a graph $G$ is denoted by $G^2$. A subset of $S\subseteq V$ is said to be an \emph{independent set} of $G$ if every pair of vertices in $S$ are non-adjacent. A subset of $K\subseteq V$ is said to be a \emph{clique} of $G$ if every pair of vertices in $K$ are adjacent. The cardinality of a clique of maximum size is called \emph{clique number} of $G$, and it is denoted by $\omega(G)$. 

A graph is called \emph{bipartite} if its vertex set can be partitioned into two independent sets. A bipartite graph $G=(X\cup Y,E)$ is called a \textit{bipartite chain graph} if there exists an ordering of vertices of $X$ and $Y$, say $\sigma_X= (x_1,x_2, \ldots ,x_{n_1})$ and $\sigma_Y=(y_1,y_2, \ldots ,y_{n_2})$, such that $N(x_{n_1})\subseteq N(x_{n_1-1})\subseteq \ldots \subseteq N(x_2)\subseteq N(x_1)$ and $N(y_{n_2})\subseteq N(y_{n_2-1})\subseteq \ldots \subseteq N(y_2)\subseteq N(y_1)$. Such an ordering of $X$ and $Y$ is called a \emph{chain ordering}, and it can be computed in linear time \cite{heggernes2007linear}. A bipartite graph $G=(X\cup Y, E)$ is said to be a \emph{star-convex bipartite} graph if a star $S$ can be defined on the vertex set $X$ such that for every vertex $y$ in $Y$, the neighbourhood of $y$ in $G$ induces a subtree of $S$. Star-convex bipartite graphs are recognizable in linear time, and associated star $S$ can also be constructed in linear time \cite{sheng2012review}. A graph $G=(V, E)$ is said to be a \emph{split graph} if $V$ can be partitioned into an independent set $S$ and a clique $K$. 


\subsection{Basic properties of $d_2$-transitivity}
In this subsection, we present some basic properties of $d_2$-transitivity. First, we show the following bounds for $d_2$-transitivity.

\begin{lem} \label{upper_bound_d_2-transitivity}
	For any graph $G$, $\Delta(G)+1\leq Tr_{d_2}(G)\leq \min\{n,(\Delta(G))^2+ 1\}$, where $\Delta(G)$ is the maximum degree of $G$.
\end{lem}

\begin{proof}
	Let $x$ be a vertex of $G$ with degree $\Delta(G)$. Consider a vertex partition $\pi=\{V_1, V_2, \ldots, V_{\Delta(G)+1}\}$ such that each $V_i$ for $2\leq i\leq \Delta(G)+1$ contains exactly one vertex from $N_G[x]$ and all the other vertices are in $V_1$. Clearly, $\pi$ forms a $d_2$-transitive partition of $G$. Therefore, $Tr_{d_2}(G)\geq \Delta(G)+1$. 
	
	Let $Tr_{d_2}(G)= k$. Clearly, $Tr_{d_2}(G)\leq n$, where $n$ is the number of vertices of $G$. Let $\pi=\{V_1, V_2, \ldots, V_k\}$ be a $d_2$-transitive partition of $G$ of size $k$. Also, let $x\in V_k$ and $N_G(x)=\{x_1, x_2, \ldots, x_l\}$. First, we show that $\displaystyle{\sum_{i=1}^{l}deg(x_i)}\geq k-1$. If $l\geq k-1$, then we are done. Otherwise, let us assume that $l<k-1$. Hence, there are some sets in $\pi$ that do not contain any vertex from $N_G(x)$. Let $V_i$ be such a set in $\pi$ and $y\in V_i$ $d_2$-dominates $x$. This implies that $y$ is adjacent to some vertex of $N_G(x)$, and hence one vertex from $V_i$ contributes one to the sum $\displaystyle{\sum_{i=1}^{l}deg(x_i)}$. Also, if $V_j$ is a set in $\pi$ that contains a vertex from $N_G(x)$, then the vertex $x$ contributes one to the sum $\displaystyle{\sum_{i=1}^{l}deg(x_i)}$. In either case, we have a contribution of one to the sum corresponding to each set in $\pi$, except the last set, $V_k$. Hence, $\displaystyle{\sum_{i=1}^{l}deg(x_i)}\geq k-1$. Since the maximum degree is $\Delta(G)$, we have $(\Delta(G))^2\geq k-1$. Therefore,  $Tr_{d_2}(G)\leq \min\{n,(\Delta(G))^2+ 1\}$.
\end{proof}

%
%

The $d_2$-transitivity of paths and cycles in the following propositions is immediately due to the above bound.

\begin{pre}\label{Path_d_2}
	Let $P_n$ be a path of $n$ vertices, and then the $d_2$-transitivity of $P_n$ is given as follows:
	
	$Tr_{d_2}(P_n) = \begin{cases}
		1 & n=1  \\
		2 & n=2 \\
		3 & n=3,4 \\
		4 & n=5,6 \\
		5 & n\geq 7
	\end{cases}$
\end{pre}

%
%
%
%
%

\begin{pre}\label{Cycle_d_2}
	Let $C_n$ be a cycle of $n$ vertices, and then the $d_2$-transitivity of $C_n$ is given as follows:
	
	$Tr_{d_2}(C_n) = \begin{cases}
		3 & n=3 \\
		4 & n=4 \\
		5 & n\geq 5
	\end{cases}$
\end{pre}

Next, we characterize graphs with small $d_2$-transitivity.

\begin{lem}\label{d2trgeq2}
Let $G$ be a connected graph. Then we have the following results:
\begin{enumerate}
	\item[(a)] $Tr_{d_2}(G)=1$ if and only if $G=K_1$.
	
	\item[(b)] $Tr_{d_2}(G)=2$ if and only if $G=K_2$.
	
	\item[(c)] $Tr_{d_2}(G)=3$ if and only if $G\in \{P_3, K_3, P_4\}$.
\end{enumerate}

\end{lem}
\begin{proof} Proofs for (a) and (b) are trivial, so we omit them.
%

(c) If $G\in \{P_3, K_3, P_4\}$, then clearly, $Tr_{d_2}(G)=3$. Conversely, let $Tr_{d_2}(G)=3$. Now by Lemma \ref{upper_bound_d_2-transitivity}, we have $\Delta(G)+1\leq 3$, that is, $\Delta(G)\leq 2$. Therefore, $G$ is either a path or a cycle. From Proposition \ref{Path_d_2}, we know that $P_3$ and $P_4$ are the only two paths for which $d_2$-transitivity is $3$. On the other hand, from Proposition \ref{Cycle_d_2}, we know that $C_3$ (equivalently $K_3$) is the only cycle for which $d_2$-transitivity is $3$. Therefore, if $Tr_{d_2}(G)=3$, then $G\in \{P_3, K_3, P_4\}$.
\end{proof}

We also characterize the graphs as having $d_2$-transitivity equal to $n$, where $n$ is the number of vertices of the graph.

\begin{lem}\label{d_2-transitivity_diam2}
Let $G$ be a graph with $n$ vertices. Then $Tr_{d_2}(G)=n$ if and only if $diam(G)\leq 2$.
\end{lem}

\begin{proof}
If $Tr_{d_2}(G)=n$, then every vertex in $G$ $d_2$-dominates every other vertex in $G$. Therefore, the distance between every pair of vertices is at most two. Therefore, $diam(G)\leq 2$.

On the other hand, if $diam(G)\leq 2$, then the distance between every pair of vertices is at most two. Hence, by putting every vertex in different sets, we get a $d_2$-transitive partition of size $n$. Therefore, $Tr_{d_2}(G)=n$.
\end{proof}


\begin{remk}	
Many important graph classes like threshold graphs, $(2K_2, P_4)$-free graphs, connected strongly regular graphs, etc. have a diameter of at most two. The lemma \ref{d_2-transitivity_diam2} implies that for these graph classes, we can solve \textsc{M$d_2$TP} trivially.
\end{remk}

%
%
%
%
%
%
%

\section{Algorithms for \textsc{M$d_2$TP}}

In this section, we find the $d_2$-transitivity for the complement of bipartite graphs and bipartite chain graphs.

\subsection{The complement of bipartite graphs}

In this subsection, we find the $d_2$-transitivity of the complement of bipartite graphs. Let $G$ be the complement of a bipartite graph $\overline{G}=(X\cup Y, \overline{E})$.

\begin{lem}
Let $G$ be the complement of a bipartite graph $\overline{G}=(X\cup Y, \overline{E})$ with $|X|=n$ and $|Y|=m$. Also let $X'=\{x\in X~|~deg_{\overline{G}}(x)=m \}$ and $Y'=\{y\in Y~|~deg_{\overline{G}}(y)=n \}$. If $\overline{G}[X' \cup Y']$ has a maximum matching of size $t$, then $Tr_{d_2}(G)=n+m-t$.
%
%
\end{lem}

\begin{proof}
Let $M=\{e_1, e_2, \ldots, e_t\}$ be a maximum matching in $\overline{G}[X' \cup Y']$ of size $t$ and $e_i=x_iy_i$ for all $1\leq i\leq t$. Let $X_t=\{x_1, x_2, \ldots, x_t\}$ and $Y_t=\{y_1, y_2, \ldots, y_t\}$. Note that since $\overline{G}[X' \cup Y']$ forms a complete bipartite graph, $M$ saturates either $X'$ or $Y'$. So, without loss of generality, let us assume that $Y_t=Y'$. Consider a vertex partition, say $\pi=\{V_1, V_2, \ldots, V_{n+m-t}\}$, of $G$ of size $n+m-t$ as follows: $V_i=\{x_i, y_i\}$ for all $1\leq i\leq t$, and every $V_j$ contains exactly one vertex from $(X\setminus X_t) \cup (Y\setminus Y_t)$ for all $t+1\leq j\leq n+m-t$. We show that $\pi$ is a $d_2$-transitive partition of $G$. Note that every vertex of $Y\setminus Y_t$ is adjacent to at least one vertex of $X\setminus X'$ in $G$. Therefore, every pair of vertices of $(X\setminus X_t) \cup (Y\setminus Y_t)$ are within distance two from each other. Therefore, $V_p$ $d_2$-dominates $V_q$ for all $t+1\leq p<q \leq n+m-t$. Also, since every set in $\{V_1, V_2, \ldots, V_t\}$ contains vertices from $X$ and $Y$ both, $V_i$ $d_2$-dominates every set in $\pi$ for all $1\leq i\leq t$. Therefore, $\pi$ is a $d_2$-transitive partition of $G$. Hence, $Tr_{d_2}(G)\geq n+m-t$.

On the other hand, let us assume that $G$ has a $d_2$-transitive partition, say $\pi$, of size
more than $n+m-t$. Since there are $n+m-2t$ vertices in $(X\setminus X_t) \cup (Y\setminus Y_t)$, at most $n+m-2t$ many sets of $\pi$ contains vertices from $(X\setminus X_t) \cup (Y\setminus Y_t)$. Therefore, there are at least $t+1$ many remaining sets in $\pi$ that contain vertices from $X_t \cup Y_t$. By the pigeonhole principle, there are at least two sets in $\pi$, say $V_i$ and $V_j$, such that $V_i$ contains vertices only from $X_t$ and $V_j$ contains vertices only from $Y_t$. But, in that case, neither set $d_2$-dominates the other because the distance between any vertex of $X_t$ and any vertex of $Y_t$ is more than two. This contradicts the fact that $\pi$ is a $d_2$-transitive partition of $G$. So, we have $Tr_{d_2}(G)\leq n+m-t$.

Therefore, the $d_2$-transitivity of $G$ is $n+m-t$. 
\end{proof}

From the above lemma, we can compute the $d_2$-transitivity as follows: given the complement of a bipartite graph $G$, first, we compute the sizes of $X'$ and $Y'$ in linear time. The size  $t$  of the maximum matching $M$ is clearly the minimum of these sizes. Then $Tr_{d_2}(G)= n+m-t$. Hence, we have the following theorem.

\begin{theo}
The \textsc{Maximum $d_2$-Transitivity Problem} can be solved in linear time for the complement of bipartite graphs.
\end{theo}

\subsection{Bipartite chain graphs}

In this subsection, we find the $d_2$-transitivity of a given bipartite chain graph $G$. To find the $d_2$-transitivity of a given bipartite chain graph $G$, first we show that the $d_2$-transitive partition of a graph $G$ is the same as the transitivity of its square graph, namely $G^2$. Then we show that if $G$ is a connected bipartite chain graph, then $G^2$ is the complement of another bipartite chain graph $H$, that is, $G^2=\overline{H}$.

\begin{lem}\label{Square equal distance_2}
For any graph $G$, $Tr_{d_2}(G)=Tr(G^2)$.

\end{lem}

\begin{proof}
Let $Tr_{d_2}(G)=k$ and $\pi=\{V_1, V_2, \ldots, V_k\}$ be a $d_2$-transitive partition of $G$. Then, by the definition of $d_2$-transitive partition, $V_i$ $d_2$-dominates $V_j$ for all $1\leq i<j\leq k$. So, for all $y\in V_j$, there exists a $x\in V_i$, such that $d(x, y)\leq 2$ in $G$. From the definition of $G^2$, it follows that there is an edge $xy$ in $G^2$. This implies that $V_i$ dominates $V_j$ in $G^2$ and hence $\pi$ is a transitive partition of $G^2$ of size $k$. Therefore, $Tr_{d_2}(G)\leq Tr(G^2)$.

On the other hand, let $Tr(G^2)=k$ and $\pi=\{V_1, V_2, \ldots, V_k\}$ be a transitive partition of $G^2$. By the definition of transitive partition, $V_i$ dominates $V_j$ for all $1\leq i<j\leq k$. So, for all $y\in V_j$, there exists a $x\in V_i$, such that $xy$ is an edge in $G^2$. From the definition of $G^2$, it follows that $d(x, y)\leq 2$ in $G$. This implies that $V_i$ $d_2$-dominates $V_j$ in $G$, and hence $\pi$ is a $d_2$-transitive partition of $G$ of size $k$. Therefore, $Tr_{d_2}(G)\geq Tr(G^2)$. So, we have $Tr_{d_2}(G)=Tr(G^2)$.
\end{proof}

%

Next, we show that if $G$ is a connected bipartite chain graph, then $G^2=\overline{H}$, where $H$ is another bipartite chain graph.

\begin{lem}\label{Square_equal_complement of bipartite}
If $G$ is a connected bipartite chain graph, then $G^2$ is the complement of another bipartite chain graph, $H$.
\end{lem}

\begin{proof}
Let $G=(X\cup Y, E)$ be a connected bipartite chain graph with chain ordering $\sigma_X= (x_1,x_2, \ldots ,x_m)$ and $\sigma_Y=(y_1,y_2, \ldots ,y_n)$ such that $N(x_m)\subseteq N(x_{m-1})\subseteq \ldots \subseteq N(x_1)$ and $N(y_n)\subseteq N(y_{n-1})\subseteq \ldots \subseteq N(y_1)$. Note that every pair of vertices of $X$ are at a distance of two as $N(y_1)=X$. Hence, $X$ forms a clique in $G^2$. Similarly, $Y$ forms another clique in $G^2$. Moreover, for any two vertices $x\in X$ and $y\in Y$, the distance between $x$ and $y$ is either $1$ or greater or equal to $3$. Hence, all the edges of $G$ are present in $G^2$, and there is no other edge in $G^2$. Now, the complement of $G^2$ is a bipartite graph. Let $H$ be the bipartite graph such that $\overline{H}=G^2$. Next, we prove that $H$ is a bipartite chain graph. Since the edges across $X$ and $Y$ are the same in $G$ and $G^2$, $N(x_m)\subseteq N(x_{m-1})\subseteq \ldots \subseteq N(x_1)$ and $N(y_n)\subseteq N(y_{n-1})\subseteq \ldots \subseteq N(y_1)$ in $G^2$ as well. Therefore, in $H$, we have $N(x_1)\subseteq N(x_{2})\subseteq \ldots \subseteq  N(x_m)$ and $N(y_1)\subseteq N(y_{2})\subseteq \ldots \subseteq N(y_n)$. Therefore, vertices of $H$ have chain ordering, and hence $H$ is a bipartite chain graph.
\end{proof}

From the Lemma \ref{Square equal distance_2} and Lemma \ref{Square_equal_complement of bipartite}, we know that finding the $d_2$-transitivity of a bipartite chain graph $G$ is the same as finding the transitivity of the complement of another bipartite chain graph $H$. From the proof of Lemma \ref{Square_equal_complement of bipartite}, it follows that we can obtain $H$ by taking the complement of $G$ and then deleting edges inside $X$ and $Y$. Note that $H$ contains some isolated vertices as $N(x_m)= N(y_n)= \emptyset$. Hence, we can compute $H$ in linear time. Moreover, the transitivity of the complement of a bipartite chain graph can be computed in linear time \cite{santra2023transitivity}. Therefore, we have the following theorem:

\begin{theo}
The \textsc{Maximum $d_2$-Transitivity Problem} can be solved in linear time for bipartite chain graphs.
\end{theo}

\section{NP-completeness of \textsc{M$d_2$TDP}}
In this section, we present three NP-completeness results for \textsc{Maximum $d_2$-Transitivity Decision Problem}, namely in split graphs, bipartite graphs, and star-convex bipartite graphs.

\subsection{Split graphs}
In this subsection, we show that the \textsc{Maximum $d_2$-Transitivity Decision Problem (M$d_2$TDP)} is NP-complete for split graphs, an important subclass of chordal graphs.

\begin{theo}
The	\textsc{Maximum $d_2$-Transitivity Decision Problem} is NP-complete for split graphs.
\end{theo}

\begin{proof}
Given a vertex partition $\pi=\{V_1, V_2, \ldots, V_k\}$ of a given split graph, we can verify in polynomial time whether $\pi$ is a $d_2$-transitive partition of that graph or not. Hence, the \textsc{Maximum $d_2$-Transitivity Decision Problem (M$d_2$TDP)} is in NP. To prove that this problem is NP-hard, we show a polynomial time reduction from the \textsc{Maximum Transitivity Decision Problem} in general graphs, which is known to be NP-complete \cite{hedetniemi2004iterated}. The reduction is as follows: given an instance of the \textsc{Maximum Transitivity Decision Problem}, that is, a graph $G=(V, E)$ and an integer $k$, we first subdivide each edge of $G$ exactly once. Let $u_i$ be the subdivision vertex corresponding to the edge $e_i\in E$ for every $1\leq i\leq m$. Then we put edges between every pair of subdivision vertices. Let the new graph be $G'=(V', E')$. Clearly, $G'$ is a split graph having $n+m$ vertices and $\frac{m^2+3m}{2}$ edges. The construction of $G'$ is illustrated in Figure \ref{fig:distance_2_split_np}. Next, we prove the following claim.
\begin{figure}[htbp!]
	\centering
	\includegraphics[scale=0.85]{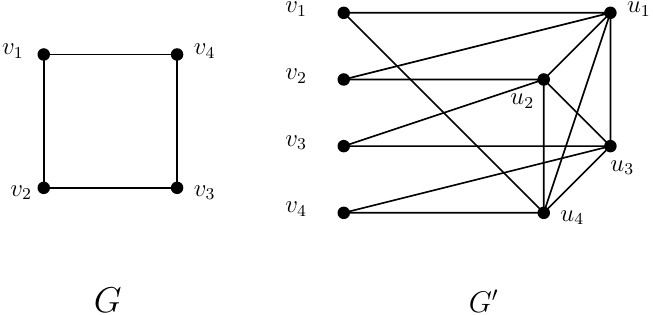}
	\caption{Construction of $G'$}
	\label{fig:distance_2_split_np}
	
\end{figure}


\begin{cl}
	The graph $G$ has a transitive partition of size $k$ if and only if $G'$ has a $d_2$-transitive partition of size $k+m$.
\end{cl}

\begin{proof}
	Let $\pi=\{V_1, V_2, \ldots, V_k\}$ be a transitive partition of $G$ of size $k$. Let us consider a vertex partition, say $\pi'=\{V'_1, V'_2, \ldots, V'_{k+m}\}$ of $G'$, as follows: $V'_i=V_i$ for all $1\leq i \leq k$ and $V'_{k+j}=\{u_j\}$, $1\leq j \leq m$. We show that $\pi'$ is a $d_2$-transitive partition of $G'$. For any pair of sets $V'_i$ and $V'_j$ with $1\leq i< j\leq k$, $V'_i$ $d_2$-dominates $V'_j$ in $G'$ as $V_i$ dominates $V_j$ in $G$. Also, the set $\{u_1, u_2, \ldots, u_m\}$ induces a complete graph in $G'$. Hence, $V'_i$ $d_2$-dominates $V'_j$ in $G'$ for $k+1\leq i< j\leq k+m$. Finally, every vertex of $G'$, other than the subdivision vertices, is adjacent to at least one subdivision vertex in $G'$. Hence, $V'_i$ $d_2$-dominates $V'_j$ in $G'$ for $1\leq i\leq k$ and  $k+1\leq j\leq k+m$. Therefore, $\pi'$ is a $d_2$-transitive partition of $G'$ of size $k+m$.

	%
		%
		%
	%

	Conversely, let $\pi=\{V_1,V_2,\ldots ,V_{k+m}\}$ be a $d_2$-transitive partition of $G'$ of size $k+m$. Let $V_{p_1}, V_{p_2}, \ldots, V_{p_t}$ be the sets in $\pi$ that do not contain any subdivision vertex, where ${p_1}< {p_2}< \ldots < {p_t}$. Since there are $m$ subdivision vertices in $G'$, there exist at least $k$ such sets in $\pi$. Therefore, $t\geq k$. Let us consider the vertex partition, say $\pi'=\{V'_1,V'_2,\ldots ,V'_{k}\}$ of $G$ as follows: $V'_i=V_{p_i}$ for $2\leq i\leq k$, and $V'_1$ contains the rest of the vertices of $G$. Since $V_{p_i}$ $d_2$-dominates $V_{p_j}$ in $G'$, every vertex of $V_{p_j}$ must be adjacent to some vertex of $V_{p_i}$ in $G$. Therefore, $\pi'$ is a transitive partition of $G$ of size $k$.
\end{proof}

%

From the above claim, it follows that the \textsc{Maximum $d_2$-Transitivity Decision Problem} is NP-complete for split graphs.
\end{proof}

\subsection{Bipartite graphs}

In this subsection, we show that the \textsc{Maximum $d_2$-Transitivity Decision Problem} is NP-complete for bipartite graphs as well.

\begin{theo}
The	\textsc{Maximum $d_2$-Transitivity Decision Problem} is NP-complete for bipartite graphs.

\end{theo}

\begin{proof}
Given a vertex partition $\pi=\{V_1, V_2, \ldots, V_k\}$ of a bipartite graph, we can verify in polynomial time whether $\pi$ is a $d_2$-transitive partition of that graph or not. Hence, the \textsc{Maximum $d_2$-Transitivity Decision Problem} is in NP. To prove that this problem is NP-hard, we show a polynomial time reduction from the \textsc{Maximum Transitivity Decision Problem} in general graphs, which is known to be NP-complete \cite{hedetniemi2004iterated}. The reduction is as follows: given an instance of the \textsc{Maximum Transitivity Decision Problem}, that is, a graph $G=(V, E)$ and an integer $k$, we construct another graph $G'=(V', E')$. Let $V=\{v_1, v_2, \ldots, v_n\}$ and $E= \{e_1, e_2, \ldots, e_m\}$ vertices and edges of $G$, respectively. For each vertex $v_i\in V$, we take two vertices $v^1_i$ and $v^2_i$ in $V'$, and for each edge $e_j\in E$, we take two vertices $u^1_j$ and $u^2_j$ in $V'$. Hence, $V'=V_1 \cup V_2\cup U_1\cup U_2$, where $V_1=\{v^1_1, v^1_2, \ldots, v^1_n\}$, $V_2=\{v^2_1, v^2_2, \ldots, v^2_n\}$, $U_1=\{u^1_1, u^1_2, \ldots, u^1_m\}$ and $U_2=\{u^2_1, u^2_2, \ldots, u^2_m\}$. Next, we add edges between $v^1_i$ and $u^1_j$ if $e_j$ is incident on $v_i$ in $G$. Similarly, we add edges between $v^2_i$ and $u^2_j$ if $e_j$ is incident on $v_i$ in $G$. Finally, we add edges between every vertex of $U_1$ and every vertex of $U_2$; that is, $U_1$ and $U_2$ induce a complete bipartite graph in $G'$. Clearly, $G'$ is a bipartite graph with $V_1\cup U_2$ and $V_2\cup U_1$ forming the bipartition, and $G'$ has $2(m+n)$ vertices and $m^2+4m$ edges. The construction is illustrated in Figure \ref{fig:distance_2_bipartite_np}. Next, we prove the following claim.

\begin{figure}[htbp!]
	\centering
	\includegraphics[scale=0.85]{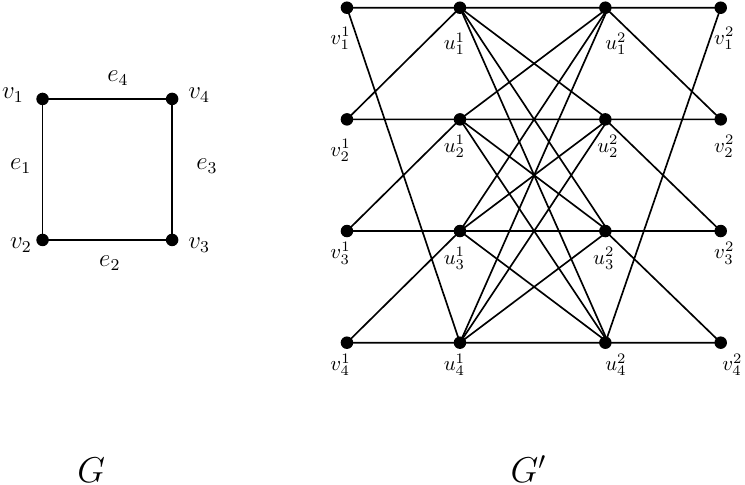}
	\caption{Construction of $G'$}
	\label{fig:distance_2_bipartite_np}
\end{figure}

\begin{cl}
	The graph $G$ has a transitive partition of size $k$ if and only if $G'$ has a $d_2$-transitive partition of size $k+2m$.
\end{cl}
\begin{proof}
	Let $\pi=\{W_1, W_2, \ldots, W_k\}$ be a transitive partition of $G$ of size $k$. Let us consider the following vertex partition, say $\pi'=\{W'_1, W'_2, \ldots, W'_{k+2m}\}$ of $G'$ as follows: for each $1\leq i\leq k$, $W'_i$ is defined as $W'_i=\{v^1_j, v^2_j| v_j\in W_i\}$ and for each $k+1\leq i\leq 2m$, $W'_i$ contains exactly one vertex from $U_1\cup U_2$. Clearly, $\pi'$ is a vertex partition of $G'$ of size $k+2m$. We show that $\pi'$ is a $d_2$-transitive partition of $G'$. Note that, by the construction of $G'$, it follows that if $e_t=v_pv_q$ is an edge in $G$, then $d(v^1_p, v^1_q)=d(v^2_p, v^2_q)=2$ in $G'$. Hence, since $W_i$ dominates $W_j$ in $G$, $W'_i$ $d_2$-dominates $W'_j$ in $G'$ for all $1\leq i< j\leq k$. Also, since $U_1\cup U_2$ induces a complete bipartite graph in $G'$, $W'_i$ $d_2$-dominates $W'_j$ for all $k+1\leq i< j\leq k+2m$. Note that every vertex of $U_1$ is at a distance of two from every vertex of $V_2$, and every vertex of $U_2$ is at a distance of two from every vertex of $V_1$. Therefore, every $W'_j$ for $k+1\leq j\leq k+2m$, is $d_2$-dominated by every $W'_i$, where $1\leq i \leq k$. Therefore, $\pi'$ is a $d_2$-transitive partition of $G'$ of size $k+2m$.

Conversely, let $\pi=\{W_1,W_2,\ldots ,W_{k+2m}\}$ be a $d_2$-transitive partition of $G'$ of size $k+2m$. Let $W_{p_1}, W_{p_2}, \ldots, W_{p_t}$ be the sets in $\pi$ that do not contain any vertex from $U_1\cup U_2$, where ${p_1}<{p_2}< \ldots< {p_t}$. Since there are $2m$ vertices in $U_1\cup U_2$, there exist at least $k$ such sets in $\pi$. Therefore, $t\geq k$. Note that since the distance between any vertex of $V_1$ and any vertex of $V_2$ is more than two, if $W_{p_t}$ contains a vertex from $V_1$ (or $V_2$), then $W_{p_i}$ for all $i<t$ contains at least one vertex from $V_1$ (respectively from $V_2$). Let us assume that $W_{p_t}$ contains vertices from $V_1$. Consider a vertex partition, say $\pi'=\{W'_1,W'_2,\ldots ,W'_{k}\}$, of $G$ as follows: $W'_i=\{v_j| v^1_j\in W_{p_i}\}$ for each $2\leq i\leq k$, and $W'_1$ contains every other vertex of $G$. Note that if $v^1_r\in W_{p_j}$ for some $2\leq j\leq k$, then every $W_{p_i}$, with $i<j$, contains at least one vertex, say $v^1_s$, such that $d(v^1_s, v^1_r)=2$ in $G'$. By the construction of $G'$, it follows that $v_s$ and $v_r$ are adjacent in $G$. Hence, $W'_i$ dominates $W'_j$ for all $1\leq i< j\leq k$. Therefore, $\pi'$ is a transitive partition of $G$ of size $k$. For the case where $W_{p_t}$ does not contain any vertex from $V_1$, that is, it contains vertices from $V_2$ only, we can argue in a similar way by considering vertices from $V_2$ and show that $\pi'$ is a transitive partition of $G$ of size $k$.
\end{proof} 

%

The above claim shows that the \textsc{Maximum $d_2$-Transitivity Decision Problem} is NP-complete for bipartite graphs.
\end{proof}

\subsection{Star-convex bipartite graphs}
In this subsection, we strengthen the NP-completeness result by showing that the \textsc{Maximum $d_2$-Transitivity Decision Problem} remains NP-complete for star-convex bipartite graphs.

\begin{theo}
The	\textsc{Maximum $d_2$-Transitivity Decision Problem} is NP-complete for star-convex bipartite graphs.
\end{theo}
\begin{proof}
Given a vertex partition $\pi=\{V_1, V_2, \ldots, V_k\}$ of a star-convex bipartite graph, we can check in polynomial time whether $\pi$ is a $d_2$-transitive partition of that graph or not. Hence, the \textsc{Maximum $d_2$-Transitivity Decision Problem} is in NP. To prove that this problem is NP-hard, we show a polynomial time reduction from the \textsc{Maximum $d_2$-Transitivity Decision Problem} in bipartite graphs. The reduction is as follows: given a bipartite graph $G=(X\cup Y, E)$ and an integer $k$, we add a vertex $x$ to $X\cup Y$ and join $x$ to every vertex of $Y$ to construct the new graph $G'=(X'\cup Y', E')$, where $X'=X \cup \{x\}$ and $Y'= Y$ and $E'= E \cup \{xy| y\in Y\}$. We can define a star $S$ on the set $X'$ with $x$ as the centre vertex. Now, if we take any vertex $y\in Y'$, $N_{G'}(y)$ induces a substar of $S$. Hence, the constructed graph $G'$ is a star-convex bipartite graph. The construction is illustrated in Figure \ref{fig:starconvexnp}. Next, we prove the following claim.

\begin{figure}[htbp!]
\centering
\includegraphics[scale=0.70]{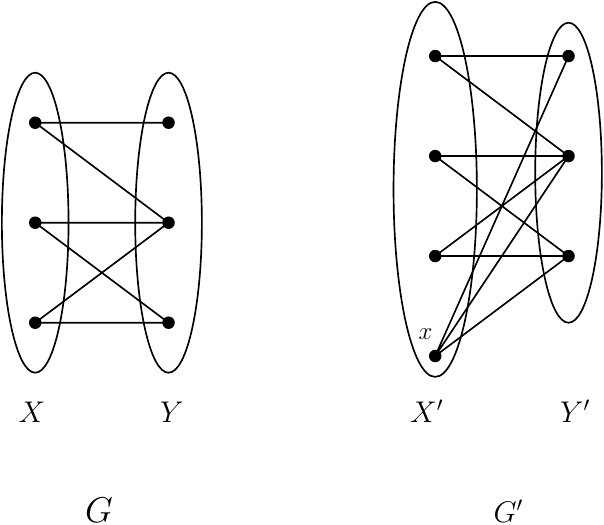}
\caption{Construction of $G'$}
\label{fig:starconvexnp}

\end{figure}

\begin{cl}
The graph $G$ has a $d_2$-transitive partition of size $k$ if and only if $G'$ has a $d_2$-transitive partition of size $k+1$.
\end{cl}
\begin{proof}
Let $\{V_1, V_2, \ldots, V_k\}$ be a $d_2$-transitive partition of $G$ of size $k$. Then, clearly, $\{V_1, V_2, \ldots, V_k, \{x\}\}$ is a $d_2$-transitive partition of $G'$ of size $k+1$ as $\{x\}$ is $d_2$-dominated by every $V_i$ in $G'$.

Conversely, let $\pi'=\{V'_1, V'_2, \ldots, V'_{k+1}\}$ be a $d_2$-transitive partition of $G'$ of size $k+1$. If $x\in V'_1$, then consider the partition $\pi=\{V_1, V_2, \ldots, V_k\}$ of $G$, where $V_1=V'_2\cup (V'_1\setminus \{x\})$ and $V_i=V'_{i+1}$ for all $2\leq i\leq k$. Since $V'_2$ $d_2$-dominates $V_j$ for all $j\geq 3$ in $G'$, $V_1$ also $d_2$-dominates every set of $\pi$ in $G$. Hence, $\pi$ is a $d_2$-transitive partition of $G$ of size $k$. On the other hand, if $x\in V'_i$ where $i\neq 1$, then consider the partition $\pi=\{V_1, V_2, \ldots, V_k\}$ of $G$, where $V_1=V'_1\cup (V'_i\setminus \{x\})$, $V_p=V'_p$, $2\leq p\leq i-1$ and $V_q=V'_{q+1}$, $i\leq q\leq k$. Since $V'_1$ $d_2$-dominates every set of $\pi'$ in $G'$, $V_1$ also $d_2$-dominates every set of $\pi$ in $G$. Hence, $\pi$ is a $d_2$-transitive partition of $G$ of size $k$.
\end{proof}

The above claim shows that \textsc{Maximum $d_2$-Transitivity Decision Problem} is NP-complete for star-convex bipartite graphs.
\end{proof}\

\section{Conclusion}
In this paper, we have introduced the notion of $d_2$-transitivity in graphs, which is a generalization of transitivity. First, we have shown some basic properties for $d_2$-transitivity. We have also proved the $d_2$-transitivity can be computed in linear time for the complement of bipartite graphs and bipartite chain graphs. On the negative side, we have shown that the \textsc{Maximum $d_2$-Transitivity Decision Problem} is NP-complete for split graphs, bipartite graphs, and star-convex bipartite graphs. It would be interesting to investigate the complexity status of this problem in other graph classes. Designing an approximation algorithm for this problem would be another challenging open problem.

\section*{Acknowledgements:} Subhabrata Paul was supported by the SERB MATRICS Research Grant (No. MTR/2019/000528). The work of Kamal Santra is supported by the Department of Science and Technology (DST) (INSPIRE Fellowship, Ref No: DST/INSPIRE/ 03/2016/000291), Govt. of India.

\bibliographystyle{alpha}
\bibliography{Dist_2_Bib}

\end{document}